

\baselineskip=14pt
\parskip=10pt
\def\halmos{\hbox{\vrule height0.15cm width0.01cm\vbox{\hrule height
  0.01cm width0.2cm \vskip0.15cm \hrule height 0.01cm width0.2cm}\vrule
  height0.15cm width 0.01cm}}
\font\eightrm=cmr8 
\font\eighttt=cmtt8
\magnification=\magstephalf
\def\P{{\cal P}}
\def\H{{\cal H}}

\def\X{{\cal X}}
\def\1{{\overline{1}}}
\def\2{{\overline{2}}}
\parindent=0pt
\overfullrule=0in

\def\frac#1#2{{#1 \over #2}}
\bf
\centerline
{
The Amazing $3^n$ Theorem  and its even more Amazing Proof 
}
\centerline
{
[Discovered by Xavier Viennot and his \'Ecole Bordelaise gang]
}
\rm
\bigskip
\centerline{ {\it
Doron 
ZEILBERGER}\footnote{$^1$}
{\eightrm  \raggedright
Department of Mathematics, Rutgers University (New Brunswick),
Hill Center-Busch Campus, 110 Frelinghuysen Rd., Piscataway,
NJ 08854-8019, USA.
{\eighttt zeilberg  at math dot rutgers dot edu} ,
\hfill \break
{\eighttt http://www.math.rutgers.edu/\~{}zeilberg/} .
Aug. 10, 2012.
Supported in part by the NSF.
Exclusively published in the Personal Journal of Shalosh B. Ekhad and Doron Zeilberger 
({\eighttt http://www.math.rutgers.edu/\~{}zeilberg/pj.html}), and {\eighttt http://arxiv.org } .
}
}

{\it  \quad\quad\quad\quad\quad\quad\quad\quad\quad\quad\quad\quad\quad\quad\quad\quad\quad
Pour mon Cher ami, Guru Xavier G. VIENNOT
}

{\bf Definition}: A {\it xavier} is a 2D tower of domino pieces where all the domino pieces at the bottom floor
are contiguous (i.e. no gaps), and every domino-piece at a higher floor is placed in such a way that
its middle-line is aligned with the left-end and/or right-end of (one or two) pieces at the floor right below it.

{\bf Theorem} [Gouyou-Beauchamps \& Viennot [GV], announced in the {\it historic} article [V], p.233]: \hfill\break
The number of xaviers with $n+1$ pieces is $3^n$.

{\bf Proof} [B\'etr\'ema \& Penaud[BeP], briefly sketched (``pictorially'') on pp. 47-48 of Mireille Bousquet-M\'elou's
great {\it habilitation} thesis[Bo1], and section 3.4.2 in yet-another {\it historic} article [Bo2],
and also reproduced on p. 81 of the Flajolet-Sedgewick [FS] {\it bible}].

Define the weight of a xavier, $x$,  by: $weight(x):=z^{NumberOfPiecesOf \,\, x}$ \quad .

Let a {\it pyramid} be a xavier whose bottom floor only has one piece, and let a {\it half-pyramid} be a pyramid
where no piece is strictly to the left of the bottom piece.

Every half-pyramid is either

(i) the singleton xavier (let's call it $Z$)

(ii) a half-pyramid where all the pieces of the higher floors are {\bf strictly} to the {\bf right} of the bottom piece

(iii) a half-pyramid where at least one of the  leftmost pieces of a non-bottom floor is aligned with the left-end of the bottom piece.

Let's call the set of half-pyramids $\H$.

The set of members of (ii) are in bijection with the Cartesian product $\{Z\} \times \H$ (just remove the bottom piece,
and to go back, move the half-pyramid one unit to the right and stick $Z$ at the bottom.

Let $h$ be a half-pyramid of case (iii). We shall map it to a  triple $[Z,h_1,h_2]$
where $h_1$ and $h_2$ are smaller half-pyramids.

Let $Z_2$ be the lowest piece strictly above the bottom floor piece (that we call $Z$) that
is aligned with it. If there is nothing above $Z_2$, 
let $h_1$ be the half-pyramid obtained from $h$ by removing both
$Z$ and $Z_2$, and let $h_2=Z_2$. Otherwise, look at the lowest piece ``above'' $Z_2$ (either touching it
or not), and then, keep adding to  $h_2$ by looking at
the pieces above this newly-acquired piece, and continue recursively,
for each newly-acquired piece, finding those pieces above them (if they exist). Continue until you
can't find anything more. $h_1$ is the half-pyramid obtained by removing the bottom $Z$
and all the  pieces of $h_2$, and (the final version of) $h_2$ consists of the above pieces, at the same horizontal locations,
but closing the vertical gaps (by dropping pieces), if necessary.

So we have a mapping
$$
\H \rightarrow \{Z\} \,\, \bigcup \,\,  \{Z\} \times \H \,\,\bigcup\,\, \{Z\} \times \H \times \H \quad .
$$
It is easy to see that it is a bijection. For a member $[Z,h_1,h_2]$ stick $Z$  below
the bottom of $h_1$, one unit to the left, and then place $h_2$ completely above $h_1$, aligning the left of $h_2$ with
$Z$. Starting with  the bottom piece of $h_2$, and working your way up, and then left-to-right, one-by-one,
construct $h$, dynamically, by
``dropping'' each piece of $h_2$ until it ``lends'' on the current $h$. The final outcome is $h$.

So we can write
$$
\H \, \equiv \, \{Z\} \,\,\bigcup\,\,\{Z\} \times \H \,\,\bigcup\,\, \{Z\} \times \H \times \H \quad  .
$$
Let $H=H(z)$ be the sum of the weights of the members of $\H$. By taking weights, we get the {\it quadratic} equation
$$
H=z+zH+zH^2 \quad ,
\eqno(Jean)
$$
and if you know {\it advanced}, middle-school, mathematics (or were a Babylonian scholar 3000 years ago),
and know how to solve, via radicals, a quadratic equation, you would be able to get an ``explicit''
expression for $H$, and by {\tt taylor}ing, you could get the first one hundred (or whatever) coefficients, and
by {\it peeking} at Sloane, you would  find  out that the enumerating sequence is the sequence of
{\it Motzkin numbers}. But this is besides (our!) point. We want to keep everything elementary
(no square-roots, (and definitely no tayloring) please!)

What about pyramids? Let $\P$ be the set of pyramids.

Let $p$ be a pyramid.  Of course it may be a half-pyramid, but if not
we will map it into a pair $[h_1,p_1]$ where $h_1$ is a half-pyramid and $p_1$ is
a smaller pyramid.

If $p$ is not already a half-pyramid, there must be, looking from bottom-to-top, a {\it lowest} piece whose left-end is {\it strictly} to
the left (of course, one unit exactly) of the left-end of the bottom piece of $p$. Let's call  it
$Z_1$. If there is nothing above $Z_1$, that's our $p_1$. Otherwise, as above, find all the pieces above it,
and having found them, those above the latter (once again, either touching or not), and continue
recursively until there is nothing more. After removing all these pieces we are left with
a half-pyramid, let's call it $h_1$, and the removed pieces, after closing, if necessary, vertical gaps,
is a brand-new pyramid, let's call it $p_1$. Once again it is obvious that you can get $p$ back from $[h_1,p_1]$, by
shifting $p_1$ one unit to the left of the bottom-piece of $h_1$ and ``dropping'' the pieces, one-by-one,
from bottom to top, and from left-to-right.
So we have
$$
\P \, \equiv \,  \H \,\,\bigcup\,\, \H \times \P \quad .
$$

Taking weights, and letting $P=P(z)$ be the sum of the weights of all members of $\P$, we have
the {\it linear} equation:
$$
P = H +H P \quad ,
\eqno(Jean-Guy)
$$
that implies $P=H/(1-H)$, and if you were stupid enough to ``solve'' for $H$ above, 
and plug-it-in into $P=H/(1-H)$, you would 
get an even uglier expression for $P$, in terms of the same annoying radical sign (reproduced in [FS], eq. (86),
but I don't like radicals, so I will spare you!)

{\it Finally}, let $\X$ be the set of xaviers, and let $x$ be a typical xavier. If the bottom floor
only has one piece, then it is a pyramid. Otherwise the bottom consists of at least
two adjacent pieces. 
Let's ignore, for now, the rightmost piece of the
bottom floor, and look at the remaining pieces. These will form the bottom floor of a new, smaller,
xavier, let's call it $x_1$. To get the rest of $x_1$,
once again, let's look at all the pieces above them, and those above the latter, and continue
recursively. It is obvious that when you remove all these pieces from $x$ you would get
a half-pyramid (whose bottom was the rightmost piece of $x$), and $x_1$ is formed by closing all
the vertical gaps, getting a brand-new xavier, whose bottom floor has one-less-piece than the
bottom floor of $x$. This gives the weight-preserving bijection
$$
\X \, \leftrightarrow \,  \P \,\,\bigcup\,\, \H \times \X \quad ,
$$
and taking weights, letting $X=X(z)$ be the sum of the weights of the members of $\X$, we get
yet another {\it linear} equation
$$
X = P +H X \quad \quad ,
\eqno(Xavier)
$$
that implies that $X=P/(1-H)$, that, in turn, using the previously established $P=H/(1-H)$, imples that
$$
X=\frac{H}{(1-H)^2} \quad .
$$
Now, if you are a high-school-algebra whiz, you can take the ``explicit'' expression for
$H$ alluded to above (but intentionally suppressed here) and ``do the algebra'', and
{\it mirabile dictu}, you would get
$$
X(z)=\frac{z}{1-3z} \quad.
$$
That, in turn, implies, by taking the coefficient of $z^{n+1}$ in $X(z)$, that the number of xaviers with $n+1$ pieces
is $3^{n}$. But this is way too advanced! Here is an alternative, elementary, proof, not using the advanced
quadratic formula and tedious manipulations with radicals. We first multiply top and bottom by $z$, getting
$$
X=\frac{zH}{z(1-H)^2} \quad .
$$
Eq. $(Jean)$ implies $zH^2=H-z-zH$, and using $(1-H)^2=1-2H+H^2$, we get:
$$
X=\frac{zH}{z(1-H)^2}=\frac{zH}{z-2zH+zH^2} =\frac{zH}{z-2zH+(H-z-zH)} =
\frac{zH}{H-3zH} =\frac{zH}{(1-3z)H} = \frac{z}{1-3z} \quad ,
$$
{\it et voil\`a}, a {\it truly} elementary proof! Or is it? Not really! Since it uses the very advanced, and
counter-intuitive, notion of {\it subtraction}. Here is an even better proof, only using addition,
of the equivalent fact
$$
X=z+3zX \quad .
$$
Indeed, thanks to $(Jean)$:
$$
X=\frac{H}{(1-H)^2}= \frac{z+zH+zH^2}{(1-H)^2} 
=\frac{z}{(1-H)^2} +z\frac{H}{(1-H)^2} + \frac{zH^2}{(1-H)^2}
=\frac{z}{(1-H)^2} +zX + \frac{z H^2}{(1-H)^2}
\quad.
$$
Now using the deep identity $H=1+\frac{H}{1-H}$ twice yields
$$
\frac{1}{(1-H)^2}=1+\frac{H}{1-H}+ \frac{H}{(1-H)^2} \quad .
$$
So
$$
X=\frac{z}{(1-H)^2} +zX + \frac{z H^2}{(1-H)^2}=
z+ z \frac{H}{1-H} + z \frac{H}{(1-H)^2}  + zX + \frac{zH^2}{(1-H)^2} =
z+ z \frac{H}{1-H} + zX  + zX + z\frac{H^2}{(1-H)^2} 
$$
$$
=z+2zX+ z \left ( \frac{H}{1-H} + \frac{H^2}{(1-H)^2} \right ) =
z+2zX+  z \frac{H}{(1-H)^2} =z+2zX +zX= z+3zX \quad  \halmos \quad .
$$
I claimed above that the above algebraic proof is {\it minus-sign-free}, but one may object that in the
{\it denominators} we have $1-H$ and $(1-H)^2$ all over the place, {\bf but}
$\frac{1}{1-H}$ is just shorthand for the minus-sign-free geometric series $\sum_{i=0}^{\infty} H^i$ and
$\frac{H}{1-H}$ is $\sum_{i=1}^{\infty} H^i$ and $\frac{H}{(1-H)^2}$ is short for
$\left ( \sum_{i_1=0}^{\infty} H^{i_1} \right ) H  \left ( \sum_{i_2=0}^{\infty} H^{i_2} \right )$.

{\bf A Natural bijection between xaviers with n+1 pieces and words of length n in the alphabet} $\{-1,0,1\}$

The {\it first} proof of the $3^n$ theorem, in [GV], was via such a bijection. Alas, it was rather complicated and {\it ad hoc}.
By taking the B\'etr\'ema \& Penaud constructions above, and by reverse-engineering the above simple 
algebraic proof, one can easily get a bijection, that may or may not be the same, or equivalent to, the original bijection,
but is much easier to state.

The first step is to note that by iterating the $\P \equiv \H \,\,\bigcup\,\, \H \times \P$ bijection, there is
a natural bijection between $\P$ and non-empty lists  of half-pyramids.
By iterating $\X \equiv \P \,\,\bigcup\,\, \H \times \X$ we get that any xavier corresponds uniquely to a creature of the form
$$
[ [h_1, \dots, h_a], h_0, [f_1, \dots, f_b] ] \quad ,
$$
where $a, b \geq  0$ and $h_1, \dots, h_a, h_0, f_1, \dots , f_b$ are all half-pyramids.
From this we can extract a letter of the alphabet $\{-1,0,1\}$ and such a creature with one
less piece. Here goes.

If $a=b=0$  and $h_0$ is the single-piece $Z$, then it is rock-bottom, and the output is the empty-word.

Apply the first  B\'etr\'ema-Penaud mapping (let's call it $BP_1$) to $h_0$. If $BP_1(h_0)$ is in $\{Z\} \times \H$, i.e. of the
form $[Z,h'_0]$ for some half-pyramid $h'_0$, then the output is the letter $0$ and the smaller creature is:
$$
[ [h_1, \dots, h_a], h'_0, [f_1, \dots, f_b] ] \quad .
$$

If $BP_1(h_0)$ is in $\{Z\} \times \H \times \H$, i.e. of the form $[Z,h'_0,h''_0]$ for some half-pyramids
$h'_0,h''_0$, then the output is the letter $1$ and the smaller creature is:
$$
[ [h_1, \dots, h_a, h'_0], h''_0, [f_1, \dots, f_b] ] \quad .
$$

If $BP_1(h_0)$ is the one-piece half-pyramid $Z$, then:

if $a>0$ the output is the letter $-1$ and the
smaller creature is:
$$
[ [h_1, \dots, h_{a-1}], h_a, [f_1, \dots, f_b] ] \quad ,
$$

while if $a=0$, then we must have $b>0$ and the output is the letter $1$ (again) and the
smaller creature
$$
[ [], f_1, [f_2, \dots, f_b] ] \quad .
$$

By iterating this mapping we get a natural (recursive) mapping from our creatures (that are in natural
bijection with xaviers) and words in the alphabet $\{-1,0,1\}$.
The readers are welcome to formulate the inverse mapping.

As with the original bijection in [GV], once again we have a fast way to {\it randomly generate}
a given xavier with $n+1$ pieces. Just generate a random string of length $n$ in $\{-1,0,1\}$, and
apply the reverse of the above mapping. It is easy to see that pyramids correspond to words
all whose partial sums are non-negative, and the half-pyramids to such words for which, in addition, the total sum is exactly zero.

{\bf The Maple package BORDELAISE}

Everything here, and much more, is implemented in the Maple package {\tt BORDELAISE}, downloadable directly from

{\tt http://www.math.rutgers.edu/\~{}zeilberg/tokhniot/BORDELAISE} \quad .

The front of this article,

{\tt http://www.math.rutgers.edu/\~{}zeilberg/mamarim/mamarimhtml/bordelaise.html}  \quad ,

contains  sample input and output files, some with nice pictures.

{\bf References}

[BeP] J. B\'etrema et J.-G. Penaud, {\it Animaux et arbres guingois},
Theoretical Computer Science {\bf 117} (1993),  67-89.

[Bo1] Mireille Bousquet-M\'elou, {\it Habilitation }, 1996, \hfill \break
{\tt http://www.labri.fr/Perso/\~{}bousquet/Articles/Habilitation/habilitation.ps.gz}

[Bo2] Mireille Bousquet-M\'elou, 
{\it Rational and algebraic series in combinatorial enumeration}, 
invited paper for  the International Congress of Mathematicians 2006. 
Proceedings of the ICM. Session lectures, pp. 789--826. \hfill \break
{\tt http://fr.arxiv.org/abs/0805.0588}.

[FS] P. Flajolet and R. Sedgewick, {\it ``Analytic Combinatorics''}, Cambridge University Press, 2009. \hfill \break
[freely(!) available on-line from {\tt http://algo.inria.fr/flajolet/Publications/book.pdf} ]

[GV] D. Gouyou-Beauchamps and G. Viennot, {\it Equivalence of the two dimensional directed animals
problem to a one-dimensional path problem}, Adv. in Appl. Math. {\bf 9} (1988), 334-357.

[V] G. X. Viennot, {\it Probl\`emes combinatoires pos\'es par la physique
statistique}, S\'eminaire N. Bourbaki, expos\'e $n^{o}$ 626, $36^{e}$ ann\'ee,
in Ast\'erisque $n^{o}$ 121-122 (1985) 225-246, SMF.
[available on-line]
\end